\documentclass{amsart}
\newtheorem{Theo}{Theorem}

\begin{document}
\title{The equation $\omega(n)=\omega(n+1)$}
\author[J.-C. Schlage-Puchta]{Jan-Christoph Schlage-Puchta}
\begin{abstract}
We prove that there are infinitely many integers $n$ such that $n$ and
$n+1$ have the same number of distinct prime divisors.
\end{abstract}
\maketitle
In this note we prove the following theorem:
\begin{Theo}
There are infinitely many $n$ such that $n$ and $n+1$ have the same
number of distinct prime factors.
\end{Theo} 
Our argument will be similar to the one given by
D. R. Heath-Brown\cite{HB1} for consecutive integers with the same
number of divisors. The method given there can
be applied to a variety of similar problems involving multiplicative
functions $f$ such that $f(p^e)$ is a nonconstant function of $e$ for
$e\geq 1$, but fails here, since $2^{\omega(p^e)}$ is clearly
independent of $e$. Thus different from \cite{HB1}, where arbitrary
large special sets were constructed in a systematic way, we have to
construct our special set 
numerically. Since we weren't able to find one with more than 5
elements, we have to use rather sharp sieve estimates, provided by the
following estimate, which is a special case of a theorem of
D. R. Heath-Brown\cite{HB2}: 
\begin{Theo}
Let $c_1, \ldots c_5$ be natural numbers. Then there are infinitely
many natural numbers $n$, such that 
\[
\sum_{i=1}^5 2^{\omega(c_i n+1)} \leq 57
\]
\end{Theo}

Now we prove Theorem 1. Assume that there are 5 integers $a_1, \ldots,
a_5$ with the properties that for any $1\leq i<j\leq 5$ we have $(a_i,
a_j)=|a_i-a_j|$, and $\omega\left(\frac{a_i}{|a_i-a_j|}\right) =
\omega\left(\frac{a_j}{|a_i-a_j|}\right)$. Define $A=a_1\cdots a_5$,
and set $c_i=a_i A$. Let $n$ be an integer as
described by Theorem 2. Then for some pair $i\neq j$ we have
$\omega(a_i A n+1)=\omega(a_j A n+1)$, for otherwise the left hand
side sum would be 
bounded below by $2^1+2^2+\ldots+2^5=62$. Then consider the integers
$\frac{a_j(a_i A n+1)}{|a_i-a_j|}$ and 
$\frac{a_i(a_j A n+1)}{|a_i-a_j|}$. They are obviously consecutive, so
it suffices to show that they have the same number of prime
divisors. Since $a_i|A$, $a_i$ and $a_j A n+1$ are coprime and we get
\begin{eqnarray*}
\omega\left(\frac{a_j(a_i A n+1)}{|a_i-a_j|}\right) & = &
\omega\left(\frac{a_i}{|a_i-a_j|}\right) + \omega(a_i A n+1)\\
& = & \omega\left(\frac{a_j}{|a_i-a_j|}\right) + \omega(a_j A n+1)\\
& = & \omega\left(\frac{a_i(a_j A n+1)}{|a_i-a_j|}\right)
\end{eqnarray*}
Hence for any $n$ given by Theorem 2 we obtain one pair of consecutive
integers with the same number of distinct prime divisors, and any pair
can come only from finitely many values of $n$, hence there are
infinitely many such pairs.

It remains to find $a_1, \ldots, a_5$ with the described properties.
The problem is that the ten equations to be checked are not
independent. To illustrate this, we show that in such a set of $a_i$
no three can be consecutive integers. Suppose to the contrary
that $a+3=a_2+1=a_3+2$. Then $(a_3,a_1)=a_3-a_1=2$ and $2|a_1$, and so
taking $N=a_1/2$ gives the three equations $\omega(2N)=\omega(2N+1)$,
$\omega(2N+1)=\omega(2N+2)$ and $\omega(N)=\omega(N+1)$. The first two
equations imply $\omega(2N)=\omega(2N+2)$. If $N$ was even,
$\omega(2N)=\omega(N)$ whereas $\omega(2N+2) = \omega(N+1)+1$, which
gives a contradiction, and in the same way the case $N$ odd can be
ruled out. 

Such obstructions can be avoided, if one chooses the $a_i$ in such a
way, that the differences $|a_i-a_j|$ are divisible by many different
prime factors, for by choosing additional congruence restrictions one
can avoid situations as above. On the other hand, the
differences should be fairly small, for otherwise the $a_i$ would
become very large, which would increase the computational effort to
check a given quintuple.

After some experimentation, the following set of quintuples seemed
promising: Define $b_1=8, b_2=9, b_3=12, b_4=34, b_5=576$,
$N=2^4\cdot 3^5\cdot 5^3\cdot 7^2\cdot 11^2\cdot 13\cdot 47^2\cdot
71^2\cdot 271$, $k=110245379356152833616$ and
consider the sequence of quintuples $(l\cdot N+k +
b_1, \ldots, l\cdot N+k + b_5)$. In fact, for
$l=1202$ this gave the following quintuple
\begin{eqnarray*}
a_1 & = & 135987650281178292389624\\
 & = &\quad 2^3\cdot 13\cdot 29\cdot 71^2\cdot 431\cdot
733\cdot 28311976573\\ 
a_2 & = & 135987650281178292389625\\
 & = &\quad 3^5\cdot 5^3\cdot 7\cdot 1481\cdot 3109\cdot
80737\cdot 1720429\\ 
a_3 & = & 135987650281178292389628\\
 & = &\quad 2^2\cdot 3\cdot 11\cdot 31^2\cdot 47^2\cdot
53\cdot 6899\cdot 1327224593\\ 
a_4 & = & 135987650281178292389650\\
 & = & \quad2\cdot 5^2\cdot 11^2\cdot 13\cdot 19\cdot
271\cdot 1145107\cdot 293245787\\ 
a_5 & = & 135987650281178292390192\\
 & = & \quad 2^4\cdot 3^4\cdot 7^2\cdot 47\cdot 71\cdot
271\cdot 2367951977749 
\end{eqnarray*}
Together with the factorization of the differences, i.e.
\[
\begin{array}{lclclcl}
a_2-a_1 & = & 1 &\quad& a_4-a_2 & = & 25 = 5^2\\
a_3-a_1 & = & 4 = 2^2 &\quad& a_5-a_2 & = & 567 = 3^4\cdot 7\\
a_4-a_1 & = & 26 = 2\cdot 13 &\quad& a_4-a_3 & = & 22 = 2\cdot 11\\
a_5-a_1 & = & 568 = 2^3\cdot 71 &\quad& a_5-a_3 & = & 564 = 2^2\cdot 3\cdot 47\\
a_3-a_2 & = & 3 &\quad& a_5-a_4 & = & 542 = 2\cdot 271\\
\end{array}
\]
it is easy to check that these values of $a_i$ satisfy all necessary
properties, and therefore Theorem 1 is proven.

The computations were performed using Mathematica 4.1.

I would like to thank D. R. Heath-Brown for making me aware of this
problem and helping with its solution.

\end{document}